\documentclass[10pt,twoside]{amsart}

\usepackage{amsthm}
\usepackage{amsmath}
\usepackage{amssymb}
\usepackage{amsfonts}
\usepackage{latexsym}
\usepackage{enumitem}
\usepackage{mathrsfs}
\usepackage[all]{xy} \SelectTips{eu}{}
\usepackage{hyperref}

\newtheorem*{intsetup*}{Setup}

\newtheorem*{intNotations*}{Notations and Conventions}

\newtheorem{intthm}{Theorem}[]

\newtheorem*{intque*}{Question}
\newtheorem*{intexa*}{Example}

\newcommand{\numberseries}{\bfseries}   
\newlength{\thmtopspace}                
\newlength{\thmbotspace}                
\newlength{\thmheadspace}               
\newlength{\thmindent}                  
\newtheoremstyle{bfupright head,slanted body}
                {\thmtopspace}{\thmbotspace}
                {\slshape}{\thmindent}{\bfseries}{.}{\thmheadspace}
                {{\numberseries \thmnumber{#2\;}}\thmnote{#3}}

\newtheoremstyle{bfupright head,upright body}
                {\thmtopspace}{\thmbotspace}
                {\upshape}{\thmindent}{\bfseries}{.}{\thmheadspace}
                {{\numberseries \thmnumber{#2\;}}\thmnote{#3}}

\newtheoremstyle{fixed bf head,slanted body}
                {\thmtopspace}{\thmbotspace}{\slshape}
                {\thmindent}{\bfseries}{.}{\thmheadspace}
                {{\numberseries \thmnumber{#2\;}}\thmname{#1}\thmnote{ (#3)}}

\newtheoremstyle{fixed bf head,upright body}
                {\thmtopspace}{\thmbotspace}{\upshape}
                {\thmindent}{\bfseries}{.}{\thmheadspace}
                {{\numberseries \thmnumber{#2\;}}\thmname{#1}\thmnote{ (#3)}}

\newtheoremstyle{numbered paragraph}
                {\thmtopspace}{\thmbotspace}{\upshape}
                {\thmindent}{\upshape}{}{\thmheadspace}
                {{\numberseries \thmnumber{#2.}}}

\theoremstyle{bfupright head,slanted body}
\newtheorem{res}{}[section]             \newtheorem*{res*}{}

\theoremstyle{bfupright head,upright body}
\newtheorem{bfhpg}[res]{}               \newtheorem*{bfhpg*}{}

\theoremstyle{fixed bf head,slanted body}
\newtheorem{thm}[res]{Theorem}          \newtheorem*{thm*}{Theorem}
\newtheorem{prp}[res]{Proposition}      \newtheorem*{prp*}{Proposition}
\newtheorem{cor}[res]{Corollary}        \newtheorem*{cor*}{Corollary}
\newtheorem{lem}[res]{Lemma}            \newtheorem*{lem*}{Lemma}
         \newtheorem*{que*}{Question}
\theoremstyle{fixed bf head,upright body}
           \newtheorem*{setup*}{Setup}
       \newtheorem*{dfn*}{Definition}
\newtheorem{rmk}[res]{Remark}           \newtheorem*{rmk*}{Remark}
           \newtheorem*{exa*}{Example}
           \newtheorem*{exer*}{Exercise}

\theoremstyle{numbered paragraph}

\newlength{\thmlistleft}        
\newlength{\thmlistright}       
\newlength{\thmlistpartopsep}   
\newlength{\thmlisttopsep}      
\newlength{\thmlistparsep}      
\newlength{\thmlistitemsep}     

\newcounter{eqc}
  {\end{list}}%

\newcounter{prt}
\newenvironment{prt}{\begin{list}{\upshape (\alph{prt})}%
    {\usecounter{prt}%
      \setlength{\leftmargin}{\thmlistleft}%
      \setlength{\labelwidth}{\thmlistleft}%
      \setlength{\rightmargin}{\thmlistright}%
      \setlength{\partopsep}{\thmlistpartopsep}%
      \setlength{\topsep}{\thmlisttopsep}%
      \setlength{\parsep}{\thmlistparsep}%
      \setlength{\itemsep}{\thmlistitemsep}}}%
  {\end{list}}%

\newcounter{rqm}
  {\end{list}}%

  {\end{list}}%

\newenvironment{prf*}[1][Proof]{%
  \begin{proof}[\bf #1]
    \setcounter{equation}{0}
    }
  {\end{proof}
}

\newcommand{\pgref}[1]{\ref{#1}}

\renewcommand{\eqref}[1]{(\pgref{eq:#1})}

\newcommand{\thmcite}[2][?]{\cite[Theorem~#1]{#2}}
\newcommand{\prpcite}[2][?]{\cite[Proposition~#1]{#2}}

\newcommand{\lemcite}[2][?]{\cite[Lemma~#1]{#2}}

\numberwithin{equation}{res}

\def\urltilda{\kern -.15em\lower .7ex\hbox{\~{}}\kern .04em}

\newcommand{\GP}{\mathsf{GP}}

\newcommand{\Proj}{\mathsf{Proj}}

\newcommand{\Mod}{\mathsf{Mod}}

\newcommand{\Ind}{\mathrm{Ind}}

\newcommand{\Rop}{R^{\sf op}}
\newcommand{\op}{\sf op}

\newcommand{\cok}{\mbox{\rm coker}}

\newcommand{\im}{\mbox{\rm Im}}

\newcommand{\Ker}[1]{\nobreak{\operatorname{Ker}#1}}

\newcommand{\Hom}{\operatorname{Hom}}

\newcommand{\is}{\cong}

   \def\soft#1{\leavevmode\setbox0=\hbox{h}\dimen7=\ht0\advance
    \dimen7 by-1ex\relax\if t#1\relax\rlap{\raise.6\dimen7
    \hbox{\kern.3ex\char'47}}#1\relax\else\if T#1\relax
    \rlap{\raise.5\dimen7\hbox{\kern1.3ex\char'47}}#1\relax
    \else\if d#1\relax\rlap{\raise.5\dimen7\hbox{\kern.9ex
    \char'47}}#1\relax\else\if D#1\relax\rlap{\raise.5\dimen7
    \hbox{\kern1.4ex\char'47}}#1\relax\else\if l#1\relax
    \rlap{\raise.5\dimen7\hbox{\kern.4ex\char'47}}#1\relax
    \else\if L#1\relax\rlap{\raise.5\dimen7\hbox{\kern.7ex
    \char'47}}#1\relax\else\message{accent \string\soft
    \space #1 not defined!}#1\relax\fi\fi\fi\fi\fi\fi}

\begin{document}

\title[Construction of Gorenstein projective modules over tensor rings]%
{Construction of Gorenstein projective modules over tensor rings}

\author[G.Q. Zhao]{Guoqiang Zhao}
\address{Guoqiang Zhao: Department of Mathmatics,
Hangzhou Dianzi University, Hangzhou 310018, China}
\email{gqzhao@hdu.edu.cn}

\author[J.X. Sun]{Juxiang Sun}
\address{Juxiang Sun (Corresponding author):
School of Mathematics and Statistics, Shangqiu Normal University, Shangqiu 476000, China}
\email{Sunjx8078@163.com}

\thanks{This research was partially supported by the National Natural Science Foundation of China (Grant No. 12061026) and Foundation for University Key Teacher by Henan Province (2019GGJS204).}





\begin{abstract}
For a tensor ring $T_R(M)$, we obtain sufficient and necessary conditions to describe 
all complete projective resolutions and all Gorenstein projective modules.
As a consequence, we provide a method for constructing Gorenstein projective modules over $T_R(M)$ from the ones of $R$.
Some applications to trivial ring extensions, Morita context rings and triangular matrix rings are given.
\end{abstract}

\maketitle

\thispagestyle{empty}

\section*{Introduction}
\label{Preliminaries}
\noindent

As a suitable generalization of finitely generated projective modules, 
the notion of the Gorenstein dimension zero over noetherian rings was introduced by Auslander and Bridger \cite{1969AB}.
Enochs and Jenda \cite{1995EnochsGP} extended it to an arbitrary ring and called it Gorenstein projective module.
Given a ring,  how to construct all Gorenstein projective modules is a fundamental problem in the study of homological algebra.
In a series of papers \cite{A2023, GFGN, 2022XI, 2020GPTRI, TriExtGPMao} and so on,
Gorenstein homological modules over triangular matrix rings,
trivial ring extensions and Morita context rings are investigated.

Let $R$ be a ring and $M$ an $R$-bimodule.
Recall that the tensor ring $T_R(M)=\bigoplus_{i=0}^\infty M^{\otimes_Ri}$,
where $M^{\otimes_R0}=R$ and
$M^{\otimes_R(i+1)}=M\otimes_R(M^{\otimes_Ri})$ for $i \geqslant 0$.
Examples of tensor rings include but are not limited to
trivial ring extensions, Morita context rings, triangular matrix rings and so on.
The classical homological properties of the tensor ring are studied
in \cite{1991Cohn, 2012Ample, 1975Roganov} and so on.

In \cite{GLS} Geiss, Leclerc and Schroer studied a certain kind of tensor ring, which is 1-Gorenstein, whose modules yield a characteristic-free categorification of the root system. Inspired by this work, Chen and Lu \cite{TRchen} studied the Gorenstein
homological properties of a tensor ring $T_R(M)$ for
an $N$-nilpotent $R$-bimodule $M$, where $R$ is a noetherian ring, and $M$ is finitely generated on both sides. 
They characterized Gorenstein projective $T_R(M)$-modules in terms of $R$-modules, 
and gave the relation between the Gorenstein dimensions of rings $R$ and $T_R(M)$, which generalized the corresponding results in \cite{GLS}.
Recently, these results have been extended to the case that
$R$ is an arbitrary associative ring and $M$ is not necessarily finitely generated $R$-bimodule \cite{DL}.

The main purpose of this paper is to describe all complete projective resolutions and
all Gorenstein projective modules over a tensor ring. 


\begin{intthm} \label{THM GP}
The sequence of projective $T_{R}(M)$-modules 
$$\cdots\to\Ind(P^{k-1})\stackrel{\alpha^{k-1}}\to\Ind(P^{k})\stackrel{\alpha^{k}}\to\Ind(P^{k+1})\to\cdots$$ 
with each $\alpha^{k}$ of the form $(\ast)$ and $\alpha^{k}_{i}\in \Hom_{R}(P^{k}, F^{i-1}(P^{k+1}))$
is a complete projective resolution if and only if, for any $k\in \mathbb{Z}$, 
the following conditions are satisfied
\begin{prt}
\item[(C1)]
$\sum_{i=1}^{j}F^{i-1}(\alpha^{k}_{j-i+1})\circ\alpha^{k-1}_{i} = 0$, for $j = 1, 2, \cdots, N+1.$
\item[(C2)]
For any $(x_{1}, x_{2}, \cdots, x_{N+1})\in \bigoplus_{i=1}^{N+1}F^{i-1}(P^{k})$ with 
$\sum_{i=1}^{j}F^{i-1}(\alpha^{k}_{j-i+1})(x_{i}) = 0$ for $j = 1, 2, \cdots, N+1$,
there is $(x_{1}^{'}, x_{2}^{'}, \cdots, x_{N+1}^{'})\in \bigoplus_{i=1}^{N+1}F^{i-1}(P^{k-1})$
such that for each $1\leqslant j\leqslant N+1,$ $x_{j}= \sum_{i=1}^{j}F^{i-1}(\alpha^{k-1}_{j-i+1})(x_{i}^{'})$.
\item[(C3)]
For $1\leqslant i\leqslant N+1$ and any projective $R$-module $P$, 
if $f_{i}\in\Hom_{R}(P^{k}, F^{i-1}(P))$ satisfies $\sum_{i=1}^{j}F^{i-1}(f_{j-i+1})\circ\alpha^{k-1}_{i} = 0$
for $j = 1, 2, \cdots, N+1$, then there is $g_{i}\in\Hom_{R}(P^{k+1}, F^{i-1}(P))$ such that for each $1\leqslant j\leqslant N+1$,
$f_{j}$ = $\sum_{i=1}^{j}F^{i-1}(g_{j-i+1})\circ\alpha^{k}_{i}$. 
\end{prt}
\end{intthm}

The following result provides a method to construct Gorenstein projective modules over $T_{R}(M)$ from the ones of $R$.

\begin{intthm}
Suppose that $M$ satisfies the compatibility conditions. 
If $X\in\GP(R)$, then $\Ind(X)$ $\in$ $\GP(T_{R}(M))$.
\end{intthm}

As applications, 
we describe all complete projective resolutions and
all Gorenstein projective modules over trivial extension of rings, Morita context rings and triangular matrix rings. We not only reobtain the earlier results in this direction,
but also get some new conclusions.

\section{Preliminaries}
\label{Preliminaries}

In this section, we fix some notation, recall relevant notions.
Throughout the paper, all rings are nonzero associative rings with identity and
all modules are unitary.
For a ring $R$, we adopt the convention that an $R$-module is a left $R$-module;
right $R$-modules are viewed as modules over the opposite ring $\Rop$.
We denote by $\Mod(R)$ the category of $R$-modules,
and by $\Proj(R)$ the subcategory of $\Mod(R)$
consisting of all projective $R$-modules. 

\begin{bfhpg}[\bf Tensor rings]\label{Tensor rings}
Let $R$ be a ring and $M$ an $R$-bimodule.
We write $M^{\otimes_R 0} = R$ and $M^{\otimes_R (i+1)} = M \otimes_R (M^{\otimes_R i})$
for $i\geqslant 0$.
For an integer $N \geqslant 0$, recall that $M$ is said to be $N$-\emph{nilpotent} if
$M^{\otimes_R(N+1)}=0$.
For an $N$-nilpotent $R$-bimodule $M$, we denote by $T_R(M)=\bigoplus_{i=0}^N M^{\otimes_Ri}$
the \emph{tensor ring} with respect to $M$. It is easy to check that $T_R(M)^{\op}\is T_{\Rop}(M)$.

It follows from \cite{TRchen} that the category $\Mod(T_R(M))$ of $T_R(M)$-modules
is equivalent to the category $\Gamma$ whose objects are the pairs $(X, u)$,
where $X \in \Mod(R)$ and $u: M\otimes_R X \to X$ is an $R$-homomorphism,
and morphisms from $(X, u)$ to $(X', u')$ are
those $R$-homomorphisms $f \in \Hom_R(X, X')$ such that
$f \circ u = u' \circ (M\otimes f)$.
Unless otherwise specified, in the paper, we always view a  $T_R(M)$-module as a pair $(X,u)$
with $X \in \Mod(R)$ and $u \in \Hom_R(M \otimes_R X, X)$.

We mention that a sequence
$$(X, u) \overset{f} \longrightarrow (X', u') \overset{g} \longrightarrow (X'', u'')$$
in $\Mod(T_R(M))$  is exact if and only if
the underlying sequence
$X \overset{f} \longrightarrow X' \overset{g} \longrightarrow X''$
is exact in $\Mod(R)$.
\end{bfhpg}

 Throughout this paper, we always let $M$ denote an $N$-nilpotent $R$-bimodule.

\begin{bfhpg}[\bf The forgetful functor and its adjoint]
\label{The forgetful functor and its}
There exists a \emph{forgetful functor}
\begin{center}
$U: \Mod(T_R(M))\to \Mod(R)$,
\end{center}
which maps a  $T_R(M)$-module $(X,u)$ to the underlying $R$-module $X$.
Recall from \lemcite[2.1]{TRchen} that $U$ admits a left adjoint
\begin{center}
$\Ind: \Mod(R) \to \Mod(T_R(M))$,
\end{center}
defined as follows:
\begin{itemize}
\item For an $R$-module $X$, define
$\Ind(X)=(\bigoplus_{i=0}^N(M^{\otimes_Ri}\otimes_RX),c_X)$,
where $c_X$ is an inclusion from
$M\otimes_R(\bigoplus_{i=0}^N(M^{\otimes_Ri}\otimes_RX))\cong
\bigoplus_{i=1}^N(M^{\otimes_Ri}\otimes_RX)$
to
$\bigoplus_{i=0}^N(M^{\otimes_Ri}\otimes_RX)$.
\item For an $R$-homomorphism $f: X \to Y$,
the $T_R(M)$-homomorphism $\Ind(f): \Ind(X) \to \Ind(Y)$ can be viewed as
a formal diagonal matrix with diagonal elements
$M^{\otimes_Ri}\otimes f$ for $0 \leqslant i\leqslant N$.
\end{itemize}
\end{bfhpg}

\begin{bfhpg}[\bf The stalk functor and its adjoint]
\label{The stalk functor and its adjoints}
There exists a \emph{stalk functor}
\begin{center}
$S: \Mod(R)\to \Mod(T_R(M))$,
\end{center}
which maps an $R$-module $X$ to the  $T_R(M)$-module $(X, 0)$.
The functor $S$ admits a left adjoint
$C: \Mod(T_R(M)) \to \Mod(R)$ defined as follows:
\begin{itemize}
\item For a  $T_R(M)$-module $(X,u)$, define $C((X,u))=\cok(u)$.
\item For a morphism $f: (X,u) \to (Y,v)$ in $\Mod(T_R(M))$, the morphism
$C(f) : \cok(u) \to \cok(v)$ is induced by the universal property of cokernels.
\end{itemize}
\end{bfhpg}

\section{Construction of Gorenstein projective $T_R(M)$-modules}
\label{Construction of Gorenstein projective T_R(M)-modules}

In this section, we determine all complete projective resolutions, and hence all Gorenstein projective modules over a tensor ring.

Recall from \cite{GHD} that an exact complex $P^\bullet$ of projective $R$-modules is called \emph{complete projective} 
if it is $\Hom_R(-, P)$-exact for every projective  $R$-module $P$.
An $R$-module $G$ is called \emph{Gorenstein projective}
if there exists a complete projective resolution $P^\bullet$ such that
$G \cong \ker(P^0 \to P^1)$.
We denote by $\GP(R)$ the subcategory of Gorenstein projective $R$-modules.

Let $T_{R}(M)$ be a tensor ring. For convenience, we use $F$ to denote the functor $M\otimes_{R}-$.
Suppose $\alpha: \Ind(X) \to \Ind(Y)$ is a morphism in $\Mod(T_{R}(M))$,
then  by definition $\alpha\in \Hom_{R}(\bigoplus_{i=1}^{N+1}F^{i-1}(X)$, $\bigoplus_{i=1}^{N+1}F^{i-1}(Y))$ and 
$\alpha \circ c_{X} = c_{Y} \circ F(\alpha)$. It is not hard to obtain that $\alpha$ is of the form
\begin{center}
$\begin{pmatrix}
\alpha_{1} & 0 & 0 & \cdots & 0 \\
\alpha_{2} & F(\alpha_{1}) & 0 & \cdots & 0 \\
\alpha_{3} & F(\alpha_{2}) & F^{2}(\alpha_{1}) & \cdots & 0 \\
\cdots & \cdots & \cdots & \cdots & \cdots \\
\alpha_{N+1} & F(\alpha_{N}) &  F^{2}(\alpha_{N-1}) & \cdots &  F^{N}(\alpha_{1}) \\
\end{pmatrix}_{(N+1)\times (N+1)} (\ast)$
\end{center}
where $\alpha_{i}\in\Hom_{R}(X, F^{i-1}(Y))$.

\begin{lem} \label{lemma1 }
The sequence $\Ind(X)\stackrel{\alpha}\to\Ind(Y)\stackrel{\beta}\to\Ind(Z)$
with $\alpha$, $\beta$ of the form $(\ast)$ 
is exact if and only if the following conditions are satisfied
\begin{prt}
\item[(C1)]
$\sum_{i=1}^{j}F^{i-1}(\beta_{j-i+1})\circ\alpha_{i} = 0$, for $j = 1, 2, \cdots, N+1.$
\item[(C2)]
If $(y_{1}, y_{2}, \cdots, y_{N+1})\in \bigoplus_{i=1}^{N+1}F^{i-1}(Y)$ with 
$\sum_{i=1}^{j}F^{i-1}(\beta_{j-i+1})(y_{i}) = 0$ for $j = 1, 2, \cdots, N+1,$
then there is $(x_{1}, x_{2}, \cdots, x_{N+1})\in \bigoplus_{i=1}^{N+1}F^{i-1}(X)$
such that for each $1\leqslant j\leqslant N+1,$ $y_{j}= \sum_{i=1}^{j}F^{i-1}(\alpha_{j-i+1})(x_{i})$.
\end{prt}
\end{lem}

\begin{prf*}
Note that 
\begin{small}
$$\beta\circ \alpha =\begin{pmatrix}
\beta_{1}\alpha_{1} & 0  & \cdots & 0 \\
\beta_{2}\alpha_{1}+F(\beta_{1})\alpha_{2} & F(\beta_{1}\alpha_{1}) & \cdots & 0 \\
\beta_{3}\alpha_{1}+F(\beta_{2})\alpha_{2}+F^{2}(\beta_{1})\alpha_{3} & F(\beta_{2}\alpha_{1}+F(\beta_{1})\alpha_{2}) & \cdots & 0 \\
\cdots & \cdots & \cdots &  \cdots \\
\sum_{i=1}^{N+1}F^{i-1}(\beta_{N-i+2})\circ\alpha_{i} & F(\sum_{i=1}^{N}F^{i-1}(\beta_{N-i+1})\circ\alpha_{i}) &   \cdots &  F^{N}(\beta_{1}\alpha_{1}) \\
\end{pmatrix}$$
\end{small}
Thus $\beta\circ \alpha =0$ if and only if the elements of first column in the matrix above are all zero, 
which is the condition (C1).

Assume $(y_{1}, y_{2}, \cdots, y_{N+1})\in$ $\Ker\beta$.
Then $\sum_{i=1}^{j}F^{i-1}(\beta_{j-i+1})(y_{i}) = 0$ for $j = 1, 2, \cdots, N+1$. 
By the condition (C2), there is $(x_{1}, x_{2}, \cdots, x_{N+1})\in \bigoplus_{i=1}^{N+1}F^{i-1}(X)$
such that $(y_{1}, y_{2}, \cdots, y_{N+1}) =$ $\alpha(x_{1}, x_{2}, \cdots, x_{N+1})$.
This means that $\Ker\beta\subseteq\im\alpha$  if and only if condition (C2) is satisfied.
Therefore, we get the assertion.
\end{prf*}

\begin{lem} \label{lemma2 }
Let $\Ind(X)\stackrel{\alpha}\to\Ind(Y)\stackrel{\beta}\to\Ind(Z)$ be an exact sequence
with $\alpha$, $\beta$ of the form $(\ast)$. 
Then for any $R$-module $W$, it is $\Hom_{T_{R}(M)}(-, \Ind(W))$-exact if and only if the following condition is satisfied
\begin{prt}
\item[(C3)]
For $1\leqslant i\leqslant N+1$ and any $f_{i}\in\Hom_{R}(Y, F^{i-1}(W))$ with $\sum_{i=1}^{j}F^{i-1}(f_{j-i+1})\circ\alpha_{i} = 0$
for $j = 1, 2, \cdots, N+1$, there is $g_{i}\in\Hom_{R}(Z, F^{i-1}(W))$ such that for each $1\leqslant j\leqslant N+1$,
$f_{j}$ = $\sum_{i=1}^{j}F^{i-1}(g_{j-i+1})\circ\beta_{i}$. 
\end{prt}
\end{lem}

\begin{prf*}
Write $( )^{\ast}$ =$\Hom_{T_{R}(M)}(-, \Ind(W))$. Consider the sequence 
$\Ind(Z)^{\ast}\stackrel{\beta^{\ast}}\to\Ind(Y)^{\ast}\stackrel{\alpha^{\ast}}\to\Ind(X)^{\ast}$.
Because for any $h\in\Ind(Z)^{\ast}$, $\alpha^{\ast}\circ\beta^{\ast}(h)$ = $h\circ\beta\circ\alpha = 0$,
it follows that $\im \beta^{\ast}\subseteq\Ker \alpha^{\ast}$.

Assume $f\in \Ker\alpha^{\ast}$, then $f$ is the form of $(\ast)$ and $f_{i}\in\Hom_{R}(Y, F^{i-1}(W))$ with $1\leqslant i\leqslant N+1$.
Since $f\circ\alpha = 0$, thus $\sum_{i=1}^{j}F^{i-1}(f_{j-i+1})\circ\alpha_{i} = 0$ for $j = 1, 2, \cdots, N+1$.
By condition (C3), there is $g\in \Ind(Z)^{\ast}$ with $g$ as the form $(\ast)$ and $g_{i}\in\Hom_{R}(Z, F^{i-1}(W))$,
such that $f = g\circ\beta$.

Hence $\im \beta^{\ast} = \Ker \alpha^{\ast}$ if and only if condition (C3) is satisfied. 
\end{prf*}

From \cite[Lemma 1.9]{DL} we know that all projective $T_R(M)$-modules are of the form $\Ind(P) = (\bigoplus_{i=1}^{N+1}F^{i-1}(P), c_P)$, 
where $P$ is a projective $R$-module. 
Now, by Lemmas 2.1 and 2.2, we can construct all complete projective resolutions over a tensor ring.

\begin{thm}
The sequence of projective $T_{R}(M)$-modules 
$$\cdots\to\Ind(P^{k-1})\stackrel{\alpha^{k-1}}\to\Ind(P^{k})\stackrel{\alpha^{k}}\to\Ind(P^{k+1})\to\cdots$$ 
with each $\alpha^{k}$ of the form $(\ast)$ and $\alpha^{k}_{i}\in \Hom_{R}(P^{k}, F^{i-1}(P^{k+1}))$
is a complete projective resolution if and only if, for any $k\in \mathbb{Z}$, 
the following conditions are satisfied
\begin{prt}
\item[(C1)]
$\sum_{i=1}^{j}F^{i-1}(\alpha^{k}_{j-i+1})\circ\alpha^{k-1}_{i} = 0$, for $j = 1, 2, \cdots, N+1.$
\item[(C2)]
For any $(x_{1}, x_{2}, \cdots, x_{N+1})\in \bigoplus_{i=1}^{N+1}F^{i-1}(P^{k})$ with 
$\sum_{i=1}^{j}F^{i-1}(\alpha^{k}_{j-i+1})(x_{i}) = 0$ for $j = 1, 2, \cdots, N+1$,
there is $(x_{1}', x_{2}', \cdots, x_{N+1}')\in \bigoplus_{i=1}^{N+1}F^{i-1}(P^{k-1})$
such that for each $1\leqslant j\leqslant N+1,$ $x_{j}= \sum_{i=1}^{j}F^{i-1}(\alpha^{k-1}_{j-i+1})(x_{i}')$.
\item[(C3)]
For $1\leqslant i\leqslant N+1$ and any projective $R$-module $P$, 
if $f_{i}\in\Hom_{R}(P^{k}, F^{i-1}(P))$ satisfies $\sum_{i=1}^{j}F^{i-1}(f_{j-i+1})\circ\alpha^{k-1}_{i} = 0$
for $j = 1, 2, \cdots, N+1$, then there is $g_{i}\in\Hom_{R}(P^{k+1}, F^{i-1}(P))$ such that for each $1\leqslant j\leqslant N+1$,
$f_{j}$ = $\sum_{i=1}^{j}F^{i-1}(g_{j-i+1})\circ\alpha^{k}_{i}$. 
\end{prt}
And any complete $T_{R}(M)$-projective resolution is of this form.
\end{thm}

What in follows is a direct consequence of Theorem 2.3, which gives the characterization of all Gorenstein projective $T_{R}(M)$-modules.

\begin{cor}
Let $T_{R}(M)$ be a tensor ring. Then 
$G\in\GP(T_{R}(M))$ if and only if 
\begin{center}
$ G = \Ker\alpha^{0} = \{(x_{1}, x_{2}, \cdots, x_{N+1})\in \bigoplus_{i=1}^{N+1}F^{i-1}(P^{0}) | 
\sum_{i=1}^{j}F^{i-1}(\alpha^{0}_{j-i+1})(x_{i}) = 0 for j = 1, 2, \cdots, N+1 \}$,
\end{center}
where $P^{k}\in\Proj(R)$ and $\alpha^{k}$ is of the form $(\ast)$ 
such that for every $k\in \mathbb{Z}$ the conditions (C1)-(C3) in Theorem 2.3 are satisfied.
\end{cor}

As a special case, we can obtain all the strongly Gorenstein projective modules
over a tensor ring, which recovers the main result of \cite{GZ} and \cite{A2024}.

\begin{cor}
Let $G$ be in $\Mod(T_{R}(M))$. Then $G$ is a strongly Gorenstein projective module if and only if 
\begin{center}
$ G = \Ker\alpha = \{(x_{1}, x_{2}, \cdots, x_{N+1})\in \bigoplus_{i=1}^{N+1}F^{i-1}(P) | 
\sum_{i=1}^{j}F^{i-1}(\alpha_{j-i+1})(x_{i}) = 0 for j = 1, 2, \cdots, N+1 \}$,
\end{center}
where $P\in\Proj(R)$, and $\alpha$ is of the form $(\ast)$ with $\alpha_{i}\in\Hom_{R}(P, F^{i-1}(P))$, 
such that the following conditions are satisfied:

\begin{prt}
\item[(SC1)]
$\sum_{i=1}^{j}F^{i-1}(\alpha_{j-i+1})\circ\alpha_{i} = 0$, for $j = 1, 2, \cdots, N+1.$
\item[(SC2)]
If $(x_{1}, x_{2}, \cdots, x_{N+1})\in \bigoplus_{i=1}^{N+1}F^{i-1}(P)$ with 
$\sum_{i=1}^{j}F^{i-1}(\alpha_{j-i+1})(x_{i}) = 0$ for $j = 1, 2, \cdots, N+1,$
then there is $(x_{1}', x_{2}', \cdots, x_{N+1}')\in \bigoplus_{i=1}^{N+1}F^{i-1}(P)$
such that for each $1\leqslant j\leqslant N+1,$ $x_{j}= \sum_{i=1}^{j}F^{i-1}(\alpha_{j-i+1})(x_{i}')$.
\item[(SC3)]
For $1\leqslant i\leqslant N+1$ and any projective $R$-module $Q$, 
if $f_{i}\in\Hom_{R}(P, F^{i-1}(Q))$ satisfies $\sum_{i=1}^{j}F^{i-1}(f_{j-i+1})\circ\alpha_{i} = 0$
for $j = 1, 2, \cdots, N+1$, then there is $g_{i}\in\Hom_{R}(P, F^{i-1}(Q))$ such that for each $1\leqslant j\leqslant N+1$,
$f_{j}$ = $\sum_{i=1}^{j}F^{i-1}(g_{j-i+1})\circ\alpha_{i}$. 
\end{prt}
\end{cor}

\begin{prf*}
Note that $G$ is strongly Gorenstein projective if and only if there is a complete projective resolution
$$\cdots\to\Ind(P)\stackrel{\alpha}\to\Ind(P)\stackrel{\alpha}\to\Ind(P)\to\cdots.$$
The assertion follows from Theorem 2.3.
\end{prf*}

It follows from \cite[Corollary 1.10]{DL} that the functor $\Ind : \Mod R\to \Mod T_{R}(M)$ preserves projective modules. 
In the following we investigate when the functors $\Ind$ preserves Gorenstein projective modules.

\begin{prp}
Let $F = M\otimes_{R}-$. If for any $1\leqslant i\leqslant N$ and $P\in\Proj(R)$, 
the functors $F^{i}$ and $\Hom_{R}(-, F^{i}(P))$ 
sends complete projective $R$-resolutions
to exact sequences of $R$-modules,
then 
$$P^{\bullet} = \cdots\to P^{-1}\stackrel{f^{-1}}\to P^{0}\stackrel{f^{0}}\to P^{1}\stackrel{f^{1}}\to\cdots$$
is a complete $R$-projective resolution if and only if
$$\Ind(P^{\bullet}) =\cdots\to\Ind(P^{-1})\stackrel{\Ind(f^{-1})}\to\Ind(P^{0})\stackrel{\Ind(f^{0})}\to\Ind(P^{1})\to\cdots$$
is a complete $T_{R}(M)$-projective resolution.
 \end{prp}

\begin{prf*}
Let $P^{\bullet} = \cdots\to P^{-1}\stackrel{f^{-1}}\to P^{0}\stackrel{f^{0}}\to P^{1}\stackrel{f^{1}}\to\cdots$
be a complete $R$-projective resolution.
From the definition of the functor $\Ind$, we know that every $\Ind(f^{k})$ is a formal diagonal matrix with diagonal
elements $F^{i}(f^{k})$ for $0\leqslant i\leqslant N$.
Since $P^{\bullet}$ and $F^{i}(P^{\bullet})$ are exact for $1\leqslant i\leqslant N$ by assumption,
one has that $\bigoplus_{i=0}^{N}F^{i}(P^{\bullet})$, and hence $\Ind(P^{\bullet})$ is exact.
For any $P\in\Proj(R)$, the adjoint pair of $(\Ind, U)$ gives rise to an isomorphism 
\begin{center}
$\Hom_{T_{R}(M)}(\Ind(P^{\bullet}), \Ind(P))$
$\cong$ $\Hom_{R}(P^{\bullet}, U(\Ind(P))) $ 
\\$= \Hom_{R}(P^{\bullet}, \bigoplus_{i=0}^{N}F^{i}(P))$
$\cong$ $\bigoplus_{i=0}^{N} \Hom_{R}(P^{\bullet}, F^{i}(P))$.
\end{center}
This implies that $\Ind(P^{\bullet})$ is a complete $T_{R}(M)$-projective resolution by assumption.

Conversely, assume that $P^{\bullet}$ is a complex such that $\Ind(P^{\bullet})$
is a complete $T_{R}(M)$-projective resolution. It follows from Theorem 2.3 that 
conditions (C1)-(C3) are satisfied.
Because every $\Ind(f^{k})$ is a formal diagonal matrix with diagonal
elements $F^{i}(f^{k})$ for $0\leqslant i\leqslant N$,
conditions (C1) and (C2) mean that $F^{i}(P^{\bullet})$ is exact for $0\leqslant i\leqslant N$.
While (C3) means that $\Hom_{R}(P^{\bullet}, F^{i}(P))$ is exact for any $0\leqslant i\leqslant N$ and projective $R$-module $P$.
Thus $P^{\bullet}$ is a complete projective resolution from the case of $i=0$.
\end{prf*}

\begin{rmk}
We refer to the above conditions as the {\it compatibility conditions} on the $R$-bimodules $M$. When $M$ is 1-nilpotent, it coincides with the notion of generalized compatible bimodule introduced in \cite[Definition 3.1]{TriExtGPMao}.
\end{rmk}

As a consequence of Proposition 2.6, we have the following result, 
which provides sufficient conditions such that the functor $\Ind$ lifts Gorenstein projective modules.

\begin{cor}
Suppose that $M$ satisfies the compatibility conditions. 
Then $\Ind(X)$ $\in$ $\GP(T_{R}(M))$ provided $X\in\GP(R)$.
\end{cor}

\section{Applications}
\label{Applications}
\noindent
In this section, we give some applications to
trivial ring extensions, Morita context rings and triangular matrix rings.

\subsection{The trivial extension of rings}
\label{The trivial extension of rings}
\noindent
Let $R$ be a ring and $M$ an $R$-bimodule.
There exists a ring $R\ltimes M$,
where the addition is componentwise
and the multiplication is given by
$(r_1, m_1)(r_2, m_2) = (r_1r_2,r_1m_2 + m_1r_2)$ for $r_1, r_2 \in R$ and $m_1, m_2 \in M$.
This ring is called the \emph{trivial extension} of
the ring $R$ by the $R$-bimodule $M$; see \cite{TRIEXT1975} and \cite{TRIEXT1971}.

Suppose that the $R$-bimodule $M$ is $1$-nilpotent,
that is, $M\otimes_RM=0$.
Then it is easy to see that the tensor ring $T_R(M)$
is nothing but the trivial ring extension $R\ltimes M$. In this case $\Ind(X) = $
$(X\oplus(M\otimes_{R}X), \begin{pmatrix}0 & 0 \\1 & 0 \\
\end{pmatrix})$.
Let $X^{k} = \Ind(P^{k})$ $= (P^{k}\oplus(M\otimes_{R}P^{k}), \begin{pmatrix}0 & 0 \\1 & 0 \\
\end{pmatrix})$, in which $P^{k}$ is a projective $R$-module.
Then $X^{k}$ is a projective $R\ltimes M$-module, and 
any projective $R\ltimes M$-module is of this form.

From Theorem 2.3 one can obtain all the complete projective resolutions and hence all Gorenstein projective modules over a trivial ring extension. This is a new conclusion.

\begin{prp} \label{GPTRM = Phi in trivial extension}
Suppose that $M$ is a $1$-nilpotent $R$-bimodule.
The sequence of $R\ltimes M$-modules 
$$\cdots\to X^{k-1}\stackrel{\alpha^{k-1}}\to X^{k}\stackrel{\alpha^{k}}\to X^{k+1}\cdots$$ 
with each $\alpha^{k}$ of the form $\begin{pmatrix} \alpha^{k}_{1} & 0 \\\alpha^{k}_{2} & M\otimes_{R} \alpha^{k}_{1}\\
\end{pmatrix}$ in which $\alpha^{k}_{1}\in \Hom_{R}(P^{k}, P^{k+1})$ and $\alpha^{k}_{2}\in \Hom_{R}(P^{k}, M\otimes_{R}P^{k+1})$
is a complete projective resolution if and only if, for any $k\in \mathbb{Z}$, 
the following conditions are satisfied
\begin{prt}
\item[(C1)]
$\alpha^{k}_{1}\circ\alpha^{k-1}_{1} = 0$ and $\alpha^{k}_{2}\circ\alpha^{k-1}_{1} + (M\otimes_{R} \alpha^{k}_{1})\circ \alpha^{k-1}_{2} = 0$.

\item[(C2)]
For any $(x_{1}, x_{2})\in P^{k}\oplus(M\otimes_{R}P^{k})$ with $\alpha^{k}_{1}(x_{1}) = 0$ and 
$\alpha^{k}_{2}(x_{1}) + (M\otimes_{R}\alpha^{k}_{1})(x_{2}) = 0$,
there is $(x_{1}^{'}, x_{2}^{'})\in P^{k-1}\oplus(M\otimes_{R}P^{k-1})$
such that $x_{1}= \alpha^{k-1}_{1}(x_{1}^{'})$ and  $x_{2}= \alpha^{k-1}_{2}(x_{1}^{'}) + (M\otimes_{R} \alpha^{k-1}_{1})(x_{2}^{'})$.

\item[(C3)]
For any projective $R$-module $P$, 
if $f_{1}\in\Hom_{R}(P^{k}, P)$ and $f_{2}\in\Hom_{R}(P^{k}, M\otimes_{R} P)$ satisfy 
$f_{1}\circ\alpha^{k-1}_{1} = 0$ and
$f_{2}\circ\alpha^{k-1}_{1} + (M\otimes_{R}f_{1})\circ\alpha^{k-1}_{2} = 0$, 
then there is $g_{1}\in\Hom_{R}(P^{k+1}, P)$ and $g_{2}\in\Hom_{R}(P^{k+1}, M\otimes_{R} P)$ such that 
$f_{1} = g_{1}\circ\alpha^{k}_{1}$ and $f_{2} = g_{2}\circ\alpha^{k}_{1} + (M\otimes_{R}g_{1})\circ\alpha^{k}_{2}$. 
\end{prt}
And any complete $R\ltimes M$-projective resolution is of this form.
\end{prp}

\subsection{Morita context rings}
\label{Morita context rings}
Let $A$ and $B$ be two rings,
and let $_AV_B$ and $_BU_A$ be two bimodules,
$\phi : U\otimes_AV \to B$ a homomorphism of $B$-bimodules,
and $\psi : V\otimes_BU \to A$ a homomorphism of $A$-bimodules.
Associated with a \emph{Morita context} $(A,B,U,V,\phi,\psi)$,
there exists a \emph{Morita context ring}
\[\Lambda_{(\phi,\psi)}=\begin{pmatrix}A & V \\ U & B\end{pmatrix}.\]
Following \thmcite[1.5]{GF1982}, one can view a $\Lambda_{(\phi,\psi)}$-module
as a quadruple $(X,Y,f,g)$ with $X \in \Mod(A)$, $Y \in \Mod(B)$,
$f \in \Hom_B(U\otimes_AX,Y)$, and $g \in \Hom_A(V\otimes_BY,X)$. 

It follows from \prpcite[2.5]{GFARTIN} that Morita context rings are trivial ring extensions
whenever both $\phi$ and $\psi$ are zero.
More precisely, consider the Morita context ring $\Lambda_{(0,0)}$.
There exists an isomorphism of rings:
$$\Lambda_{(0,0)}\overset{\cong}\longrightarrow (A\times B)\ltimes (U \oplus V)
\,\,\text{via}\,\, \begin{pmatrix}a & v \\ u & b\end{pmatrix} \mapsto
((a,b),(u,v)).$$
Thus, there exists an isomorphic functor
$$\mu: \Mod(\Lambda_{(0,0)}) \to \Mod((A\times B)\ltimes (U \oplus V))\ \mathrm{via}\
(X,Y,f,g) \mapsto ((X,Y),(g,f)),$$
where $(g,f)$ is from
$(U \oplus V)\otimes_{A\times B}(X, Y) \cong (V\otimes_BY, U\otimes_AX)$ to $(X, Y)$.

In this case for $(X, Y)\in\Mod (A\times B)$,  we  have $\Ind (X, Y) =( (X, Y)\oplus (U\oplus V)\otimes_{A\times B}(X, Y), \begin{pmatrix}0 & 0 \\1 & 0 \\
\end{pmatrix})$ $\cong$ 
$((X\oplus (V\otimes_B Y), (U\otimes_A X)\oplus Y), (\begin{pmatrix} 0 \\ 1 _{V\otimes_B Y} \\
\end{pmatrix}, \begin{pmatrix}1 _{U\otimes_A X}\\0  \end{pmatrix}))$.
Thus any projective $\Lambda_{(0,0)}$-module is of the form
$\mu^{-1}(\Ind(P, Q))$, where ${_A}P$ and ${_B}Q$ are projective modules,  i.e. $(P, Q) \in \Proj(A\times B)$.

We mention that $(U \oplus V)\otimes_{A\times B}(U \oplus V) \cong (U\otimes_AV) \oplus (V\otimes_BU)$.
Then the $A\times B$-bimodule $U \oplus V$ is $1$-nilpotent if and only if $U\otimes_AV=0=V\otimes_BU$.
Let $X^{k} = \mu^{-1}(\Ind(P^{k}, Q^{k}))$ $= (P^{k}\oplus (V\otimes_B Q^{k}), (U\otimes_A P^{k})\oplus Q^{k}, \begin{pmatrix}1 _{U\otimes_A P^{k}}\\0  \\
\end{pmatrix}, \begin{pmatrix} 0 \\ 1 _{V\otimes_B Q^{k}}) \\ \end{pmatrix})$.

Thus, we can get all the complete projective resolutions over a Morita context ring from Proposition 3.1, which is due to the main result \cite[Theorem 3.1]{A2023}.
We mention that the proof here improves the ones in \cite{A2023}.

\begin{cor} \label{con GP in morita ring}
Let $\Lambda_{(0,0)}$ be a Morita context ring with $U\otimes_AV=0=V\otimes_BU$. 
The sequence of $\Lambda_{(0,0)}$-modules 
$$\cdots\to X^{k-1}\stackrel{f^{k-1}}\to X^{k}\stackrel{f^{k}}\to X^{k+1}\to\cdots$$ 
with each $f^{k}$ of the form ($\begin{pmatrix} \tau^{k} & 0 \\ \beta^{k} & V\otimes_{B} \sigma^{k}\\
\end{pmatrix}$, $\begin{pmatrix} U\otimes_{A} \tau^{k} & \gamma^{k} \\ 0 &  \sigma^{k}\\
\end{pmatrix}$), in which $\tau^{k}\in \Hom_{A}(P^{k}, P^{k+1})$, 
$\sigma^{k}\in\Hom_{B}(Q^{k}, Q^{k+1})$ and $\beta^{k}\in$ $\Hom_{A}(P^{k}$, $V\otimes_{B}Q^{k+1})$
, $\gamma^{k}\in$ $\Hom_{B}(Q^{k}, U\otimes_{A}P^{k+1})$
is a complete projective resolution if and only if, for any $k\in \mathbb{Z}$, 
the following conditions are satisfied
\begin{prt}

\item[(C1)$^{'}$]
$\tau^{k}\circ\tau^{k-1} = 0= \sigma^{k}\circ\sigma^{k-1}$ and 
$\beta^{k}\circ\tau^{k-1}+(V\otimes_{B}\sigma^{k})\circ\beta^{k-1} = 0 =\gamma^{k}\circ\sigma^{k-1}+(U\otimes_{A}\tau^{k})\circ\gamma^{k-1}$.

\item[(C2)$^{'}$]
for any $(p, q)\in (P^{k}, Q^{k})$ and $(x, y)\in (V\otimes_{B}Q^{k}, U\otimes_{A}P^{k})$ 
with $\tau^{k}(p) = 0 = \sigma^{k}(q) $ and 
$\beta^{k}(p) + (V\otimes_{B}\sigma^{k})(x) = 0 =$ $\gamma^{k}(q) + (U\otimes_{A}\tau^{k})(y)$,
there are $(p', q')\in (P^{k-1}, Q^{k-1})$ and $(x', y')\in (V\otimes_{B}Q^{k-1}, U\otimes_{A}P^{k-1})$
such that $p =  \tau^{k}(p')$, $q = \sigma^{k}(q')$ and  
$x = \beta^{k}(p') + (V\otimes_{B}\sigma^{k})(x')$,
$y = \gamma^{k} (q') + (U\otimes_{A}\tau^{k})(y')$.

\item[(C3)$^{'}$]
For any projective $A$-module $P$ and projective $B$-module $Q$, 
and $f_{1}\in\Hom_{A}(P^{k}$, $P)$, $f_{2}\in\Hom_{B}(Q^{k}, Q)$ and 
$u_{1}\in\Hom_{A}(P^{k}$, $V\otimes_{B} Q)$, $u_{2}\in\Hom_{B}(Q^{k}, U\otimes_{A} P)$ satisfy 
$$f_{1}\circ\tau^{k-1} = 0 = f_{2}\circ\sigma^{k-1}$$ 
$$u_{1}\circ\tau^{k-1}+  (V\otimes_{B} f_{2})\circ\beta^{k} = 0 =
u_{2}\circ\sigma^{k-1}+  (U\otimes_{A}f_{1})\circ\gamma^{k},$$
then there is $g_{1}\in\Hom_{A}(P^{k+1}, P)$, 
$g_{2}\in\Hom_{B}(Q^{k+1}, Q)$ and 
$v_{1}\in\Hom_{A}(P^{k+1}$, $V\otimes_{B} Q)$,
$v_{2}\in\Hom_{B}(Q^{k+1}, U\otimes_{A}P)$such that 
$$f_{1} = g_{1}\circ\tau^{k}, f_{2} = g_{2}\circ\sigma^{k}$$  
$$u_{1} = v_{1}\circ\tau^{k} + (V\otimes_{B} g_{2})\circ\beta^{k}, u_{2} = v_{2}\circ \sigma^{k}
+ (U\otimes_{A}g_{1})\circ\gamma^{k}.$$ 

\end{prt}
And any complete $\Lambda_{(0,0)}$-projective resolution is of this form.
\end{cor}

\begin{prf*}
Because $\mu: \Mod(\Lambda_{(0,0)}) \to \Mod((A\times B)\ltimes (U \oplus V))$ is an isomorphic functor,
it is equivalent to show that the assertion holds for the sequence of $(A\times B)\ltimes (U\oplus V)$-modules 
$$\cdots\to\Ind(P^{k-1}, Q^{k-1})\stackrel{\mu(f^{k-1})}\to \mu(X^{k}) = \Ind(P^{k}, Q^{k})\stackrel{\mu(f^{k})}\to \Ind(P^{k}, Q^{k})\to\cdots  (3.1)$$ 
where $\mu(f^{k}) =$ $\begin{pmatrix} (\tau^{k}, \sigma^{k}) & 0 \\  (\beta^{k}, \gamma^{k}) & (V\otimes_{B}\sigma^{k}, U\otimes_{A}\tau^{k})\\
\end{pmatrix}$ for the sake of 
$(U \oplus V)\otimes_{A\times B}(\tau^{k}, \sigma^{k})\cong$ $(V\otimes_{B}\sigma^{k}, U\otimes_{A}\tau^{k})$.

From Proposition 3.1 we know that the sequence (3.1) is a complete projective resolution if and only if, 
for any $k\in \mathbb{Z}$, the conditions (C1)-(C3) are satisfied.
Since $\alpha^{k}_{1} = (\tau^{k}, \sigma^{k})$, $\alpha^{k}_{2} = (\beta^{k}, \gamma^{k})$,
it follows that $\alpha^{k}_{1}\circ\alpha^{k-1}_{1} = (\tau^{k}\circ\tau^{k-1}, \sigma^{k}\circ\sigma^{k-1})$  
and $\alpha^{k}_{2}\circ\alpha^{k-1}_{1} + ((U \oplus V)\otimes_{A\times B} \alpha^{k}_{1})\circ \alpha^{k-1}_{2} = $
$(\beta^{k}\circ\tau^{k-1}+(V\otimes_{B}\sigma^{k})\circ\beta^{k-1}, \gamma^{k}\circ\sigma^{k-1}+(U\otimes_{A}\tau^{k})\circ\gamma^{k-1})$.
Thus (C1) is just (C1)$^{'}$.

Since $(U \oplus V)\otimes_{A\times B}(P^{k}, Q^{k})\cong$ $(V\otimes_{B}Q^{k}, U\otimes_{A}P^{k})$ $\in\Mod(A\times B)$,
the condition (C2) is equivalent to that 
for any $(p, q)\in (P^{k}, Q^{k})$ and $(x, y)\in (V\otimes_{B}Q^{k}, U\otimes_{A}P^{k})$ with $(\tau^{k}, \sigma^{k})(p, q) = 0$ and 
$(\beta^{k}, \gamma^{k})(p, q) + (V\otimes_{B}\sigma^{k}, U\otimes_{A}\tau^{k})(x, y) = 0$,
there are $(p', q')\in (P^{k}, Q^{k})$ and $(x', y')\in (V\otimes_{B}Q^{k}, U\otimes_{A}P^{k})$
such that 
\begin{center}
$(p, q) =  (\tau^{k}, \sigma^{k})(p', q')$ and  
$(x, y) = (\beta^{k}, \gamma^{k})(p', q') + (V\otimes_{B}\sigma^{k}, U\otimes_{A}\tau^{k})(x', y')$.
\end{center}
Hence  (C2) is exactly  (C2)$^{'}$. 

The condition (C3) means that for any projective $(A\times B)$-module $(P, Q)$, 
if $(f_{1}, f_{2})\in\Hom_{A\times B}((P^{k}, Q^{k}), (P, Q))$ and 
$(u_{1}, u_{2})\in\Hom_{A\times B}((P^{k}, Q^{k}), (V\otimes_{B} Q, U\otimes_{A}P))$ satisfy 
$(f_{1}, f_{2})\circ(\tau^{k-1}, \sigma^{k-1}) = 0$ and
$(u_{1}, u_{2})\circ(\tau^{k-1}, \sigma^{k-1})+  (V\otimes_{B} f_{2}, U\otimes_{A}f_{1})\circ(\beta^{k}, \gamma^{k}) = 0$, 
then there is $(g_{1}, g_{2})\in\Hom_{A\times B}((P^{k+1}, Q^{k+1}), (P, Q))$ and 
$(v_{1}, v_{2})\in\Hom_{A\times B}((P^{k+1}, Q^{k+1})$, $(V\otimes_{B} Q$, $U\otimes_{A}P))$ such that 
\begin{center}
$(f_{1}, f_{2}) = (g_{1}, g_{2})\circ(\tau^{k}, \sigma^{k})$ and $(u_{1}, u_{2}) = (v_{1}, v_{2})\circ(\tau^{k}, \sigma^{k}) 
+ (V\otimes_{B} g_{2}, U\otimes_{A}g_{1})\circ(\beta^{k}, \gamma^{k}).$
\end{center}
That is the condition (C3)$^{'}$, as desired.
\end{prf*}

In the case of $U=0$, the Morita context ring $\Lambda_{(0,0)}$ becomes a triangular matrix ring $\Gamma=\begin{pmatrix} A & V \\ 0 & B\\
\end{pmatrix}$, and a left $\Gamma$-module is identified with a triple
$(X, Y, \varphi)$, where $X\in \Mod A$, $Y\in \Mod B$, and $\varphi: V\otimes_{B}Y \to X$ is an $A$-map.  
A $\Gamma$-map $(X, Y, \varphi) \to (X', Y', \varphi')$ is identified with a pair $(f, g)$, where $f\in\Hom_{A}(X, X')$, $g\in\Hom_{A}(Y, Y')$, such that $\varphi'(V\otimes_{B}g) = f\varphi$.

In the following, we obtain all the complete projective resolutions over a triangular matrix ring. This improves the main result of \cite{GZ},
which only gave the strongly complete projective resolutions. 

\begin{cor} \label{con GP in triangular matrix ring}
Let $\Gamma=\begin{pmatrix} A & V \\ 0 & B\\
\end{pmatrix}$
be a triangular matrix ring.
The sequence of $\Gamma$-modules 
$$\cdots\to X^{k-1}\stackrel{f^{k-1}}\to X^{k}\stackrel{f^{k}}\to X^{k+1}\to\cdots,$$ 
where $X^{k} = (P^{k}\oplus V\otimes_{B} Q^{k}, Q^{k}, \begin{pmatrix} 0\\ 1_{V\otimes_{B} Q^{k}}\\
\end{pmatrix})$ and $f^{k}$ is of the form ($\begin{pmatrix} \tau^{k} & 0 \\ \beta^{k} & V\otimes_{B} \sigma^{k}\\
\end{pmatrix}$, $\sigma^{k})$ with $\tau^{k}\in \Hom_{A}(P^{k}, P^{k+1})$, 
$\sigma^{k}\in\Hom_{B}(Q^{k}, Q^{k+1})$ and $\beta^{k}\in$ $\Hom_{A}(P^{k}$, $V\otimes_{B}Q^{k+1})$,
is a complete projective resolution if and only if, for any $k\in \mathbb{Z}$, 
the following conditions are satisfied
\begin{prt}

\item[(i)]
$\mathcal{P}^{\bullet}: \cdots\to P^{k-1}\stackrel{\tau^{k-1}}\to P^{k}\stackrel{\tau^{k}}\to P^{k+1}\to\cdots$ is a complex of projective $A$-modules, and it is  $\Hom_{A}(-, P)$-exact for any projective $A$-module $P$.

\item[(ii)]
$\mathcal{Q}^{\bullet}: \cdots\to Q^{k-1}\stackrel{\sigma^{k-1}}\to Q^{k}\stackrel{\sigma^{k}}\to Q^{k+1}\to\cdots$ is an exact complex of projective $B$-modules.

\item[(iii)]
$\beta^{k}\circ\tau^{k-1}+(V\otimes_{B}\sigma^{k})\circ\beta^{k-1} = 0$.


\item[(iv)]
If $\tau^{k}(p) = 0 $ and 
$\beta^{k}(p) + (V\otimes_{B}\sigma^{k})(x) = 0$,
then there are $(p', x')\in P^{k-1}\oplus V\otimes_{B}Q^{k-1}$ such that $p =  \tau^{k}(p')$ and  
$x = \beta^{k}(p') + (V\otimes_{B}\sigma^{k})(x')$.

\item[(v)]
For any projective $B$-module $Q$, 
and $(f, g)\in\Hom_{A}(P^{k}$, $V\otimes_{B} Q)\oplus\Hom_{B}(Q^{k}$, $Q)$ satisfy 
$g\circ\sigma^{k-1}= 0= f\circ\tau^{k-1}+  (V\otimes_{B} g)\circ\beta^{k},$
then there is  
$(f', g')\in\Hom_{A}(P^{k+1}$, $V\otimes_{B} Q)\oplus\Hom_{B}(Q^{k+1}, Q)$ such that 
$g = g'\circ\sigma^{k}$ and  
$f = f'\circ\tau^{k} + (V\otimes_{B} g')\circ\beta^{k}.$ 
\end{prt}
And any complete $\Gamma$-projective resolution is of this form.
\end{cor}

\begin{prf*}
Let $k$ be any integer. Since $U = 0$,  then $\gamma^{k} =0$, and hence the condition (C1)$^{'}$ is equivalent to that $\mathcal{P}^{\bullet}$, $\mathcal{Q}^{\bullet}$ are complexes, and (iii) holds.

In (C2)$^{'}$, the condition \lq for any $q\in Q^{k}$ with $\sigma^{k}(q) =0$,
there are $q'\in Q^{k-1}$ such that $q = \sigma^{k}(q')$\rq  is 
equal to that the complex $\mathcal{Q}^{\bullet}$ is exact. The rest condition of (C2)$^{'}$ is just condition (iv). 

Note that for any projective $A$-module $P$, $\Hom_{A}(\mathcal{P}^{\bullet}, P)$ is clearly a complex.
In (C3)$^{'}$, the condition \lq  if $f_{1}\in\Hom_{A}(P^{k}$, $P)$ such that
$f_{1}\circ\tau^{k-1} = 0$, 
then there is $g_{1}\in\Hom_{A}(P^{k+1}, P)$, 
such that $f_{1} = g_{1}\circ\tau^{k}$\rq means that the complex $\mathcal{P}^{\bullet}$ is $\Hom_{A}(-, P)$-exact.
The rest condition of (C3)$^{'}$ is condition (v).

Therefore, the conclusion holds from Corollary 3.2.
\end{prf*}

\begin{rmk}
We mention that the former results have dual versions, which means that we can describe all complete injective resolutions and Gorenstein injective modules over a tensor ring in a dual manner.
\end{rmk}



\bibliographystyle{amsplain-nodash}

\def\cprime{$'$}
  \providecommand{\arxiv}[2][AC]{\mbox{\href{http://arxiv.org/abs/#2}{\sf
  arXiv:#2 [math.#1]}}}
  \providecommand{\oldarxiv}[2][AC]{\mbox{\href{http://arxiv.org/abs/math/#2}{\sf
  arXiv:math/#2
  [math.#1]}}}\providecommand{\MR}[1]{\mbox{\href{http://www.ams.org/mathscinet-getitem?mr=#1}{#1}}}
  \renewcommand{\MR}[1]{\mbox{\href{http://www.ams.org/mathscinet-getitem?mr=#1}{#1}}}
\providecommand{\bysame}{\leavevmode\hbox to3em{\hrulefill}\thinspace}
\providecommand{\MR}{\relax\ifhmode\unskip\space\fi MR }
\providecommand{\MRhref}[2]{%
  \href{http://www.ams.org/mathscinet-getitem?mr=#1}{#2}
}
\providecommand{\href}[2]{#2}

\end{document}